\documentclass[12pt,a4paper,final]{iopart}

\usepackage{iopams}  
\usepackage{graphicx}
\usepackage{epstopdf}
\usepackage{enumitem, color, amssymb}
\usepackage[dvipsnames]{xcolor}
\usepackage{hyperref}
\usepackage{amssymb}
\usepackage{listings}
\usepackage{enumitem, color, amssymb}
\usepackage[dvipsnames]{xcolor}
\usepackage{pgfplots}
\usepackage{cancel}
\usepackage{textcomp}
\usepackage{algorithm2e}
\usepackage{subfigure}
\pgfplotsset{compat=1.12}
\usetikzlibrary{decorations.markings}

\newtheorem{theorem}{Theorem}
\newcommand{\mrd}{{\mathrm{d}}}
\newcommand{\mri}{{\mathrm{i}}}

\newtheorem{definition}{Definition}

\ifdefined\COMMENTBILL
\newcommand{\billcomment}[1]{ {\color{red} WL: #1 }}%
\else
\newcommand{\billcomment}[1]{}
\fi 

\ifdefined\COMMENTNAEEM
\newcommand{\naeemcomment}[1]{ {\color{blue} ND: #1 }}%
\else
\newcommand{\naeemcomment}[1]{}
\fi

\begin{document}
\title[Polarimetric Neutron Tomography]{Polarimetric Neutron Tomography of Magnetic Fields: Uniqueness of Solution and Reconstruction}

\author{Naeem M. Desai}
\address{School of Mathematics, University of Manchester, Manchester, UK}
\ead{naeem.desai@manchester.ac.uk}

\author{William R.B. Lionheart}
\address{School of Mathematics, University of Manchester, Manchester, UK}
\ead{william.lionheart@manchester.ac.uk}

\author{Morten Sales}
\address{DTU Physics, Technical University of Denmark, Lyngby, Denmark}
\ead{msales@fysik.dtu.dk}

\author{Markus Strobl}
\address{Paul Scherrer Institute (PSI), Villigen, Switzerland}
\ead{markus.strobl@psi.ch}

\author[cor1]{S\o{}ren Schmidt}
\address{DTU Physics, Technical University of Denmark, Lyngby, Denmark}
\ead{ssch@fysik.dtu.dk}

\begin{abstract}
We consider the problem of determination of  a magnetic field from three dimensional polarimetric neutron tomography data. We see that this is an example of a non-Abelian ray transform and that the problem has a globally unique solution for smooth magnetic fields with compact support, and a locally unique solution for less smooth fields. We derive the linearization of the problem and note that the derivative is injective. We go on to show that the linearised problem about a zero magnetic field reduces to  plane Radon transforms and suggest a modified Newton Kantarovich method (MNKM) type algorithm for the numerical solution of the non-linear problem, in which the forward problem is re-solved but the same derivative used each time. Numerical experiments demonstrate that MNKM   works for small enough fields (or large enough velocities)  and show an example where it fails to reconstruct a slice of the simulated data set. Lastly we show viewed as an optimization problem the inverse problem is  non-convex  so we expect gradient based methods may fail.
\end{abstract}
\maketitle

\section{Introduction}

Neutron Tomography is widely used in science and industry; it provides important complementary information to that given by x-rays as neutrons have zero electrical charge and can penetrate deeply into massive samples, see \cite{TropReview}. Neutrons are highly sensitive to magnetic fields owing to their magnetic moment, which makes them an indispensable tool for scattering studies to investigate magnetic phenomena in material science. At the close of the last century, it was proposed to use the magnetic moment for real space tomographic imaging investigations of magnetic structures. Simulations have been presented in \cite{FirstPaper96,AustJPhys98} of a novel method called Neutron Magnetic Tomography. 

However, a decade later an efficient, but restricted, experimental approach of Polarized Neutron Imaging was introduced to make two and three dimensional images of magnetic fields \cite{Kardjilov}. In contrast to initial proposals, the setup used does not measure the full spin rotation matrix but only a single diagonal element. Correspondingly, the method has been applied for strongly oriented fields and high symmetry cases providing significant {\em a priori} knowledge for analysis, e.g. through field modeling and simulation matching to data. In this manner magnetic fields of electromagnetic devices and electric currents, \cite{Kardjilov,JPARC,PhysicaB09,JAppPhys08,IOP2015}, but also quantum mechanical effects in superconductors \cite{Kardjilov,PhysicaB09,Dawson} could be studied successfully.

A potential important application is the study of three dimensional magnetic configurations, in particular magnetic domain structures. This is beyond the capability of previous methodologies as it requires the measurement of the full polarisation matrix. The feasibility of this was demonstrated in \cite{PhysicaB09}. In this paper we assume such a full polarimetric set-up measuring the full matrix of images, of the three perpendicular incoming polarization and analysis directions, for every projection of a full tomographic scan.

Earlier attempts at  reconstruction were presented using the {\em ad hoc} method of scalar Radon transform inversion applied to polarized neutron data \cite{FirstPaper96,AustJPhys98,Kardjilov}. 
The linearized problem about a zero magnetic field was presented in \cite{3DPNT}.
We show that the problem can be formulated as a non-Abelian ray transform and that sufficiently smooth magnetic fields are uniquely determined by polarimetric neutron tomography data. We show that the linearized problem, for small magnetic fields, has a unique solution for less regular fields and propose an iterative reconstruction algorithm which we have implemented and tested on simulated data. (This formulation as a non-Abelian ray transform and the consequent uniqueness result for non-linear and linearized problem was first presented in the conference talk \cite{AIPTalk} and in more detail in the thesis \cite{DesaiThesis}.)
This new theoretical framework lays the foundation for practical three dimensional polarimetric neutron tomography of magnetic fields (PNTMF) that will facilitate the imaging of magnetic domains in metal samples and the corresponding design of magnetic materials.

\section{Physical Framework}
\label{section:experiment_procedure}
Despite being electrically neutral, neutrons carry a magnetic moment, which is coupled to their spin vector. For an ensemble of polarized neutrons in a magnetic field it can be shown that they behave like a particle with a classical magnetic moment. Given a specific ray $ x_0 + \mathbf{v}t$ with velocity $\mathbf{v}$, we define $s=vt$, speed $v=|\mathbf{v}|\ne 0$ and direction $\xi=\mathbf{v}/v$. Parametrizing in distance we write $x(s)=x_0 +s \xi$ and with $v$ fixed a ray is parameterized by its initial point $x_0$ and direction $\xi$, as is standard in tomography. Considering a fixed ray $(x_0, \xi)$ and speed $v$ the spin vector $\sigma(s)$ satisfies
\begin{equation}
\frac{\mrd}{\mrd s} \sigma(s) = \frac{\gamma_N}{v} \sigma(s) \times B(x(s)),
\label{eqn:larmor}
\end{equation}

where $\gamma_N=-1.8324\times 10^8$ ~rad s$^{-1}$~T$^{-1}$, is the gyromagnetic ratio of the neutron and $B(x(s))$ is the magnetic field at that point on the ray

Experimentally we assume a polarized neutron beam with uniform velocity and beam size of approximately $4 \times 4 ~\rm{cm}^2$. The spatial resolution that can be achieved over such a beam is defined by the basic pinhole imaging geometry, where the spatial blur is given by $d = {lD}/{\hat{L}}$ as in \cite{TropReview}. Here $l$ is the distance from sample to detector, ${\hat{L}}/{D}$ is the collimation ratio of the distance between pinhole and sample, $\hat{L}$, and the pinhole diameter $D$. In reality the beam has a velocity spread of order $1\%$ that can be tuned by trading against the neutron flux. 

For the neutron velocity range commonly used, polarizations well beyond $90\%$ and close to $100\%$ can be achieved with spin filters such as super mirror devices as used for polarized neutron imaging \cite{Kardjilov,Dawson}. A neutron spin filter only transmits neutrons with spin parallel to the magnetization of the device. Assuming the beam is not depolarized locally within the setup, a well defined final Larmor precession angle with respect to the analyzer's direction can be extracted. 

For our simplified model,  effects such as spin precession being wavelength dependent, initial polarization and conservation in the setup without the sample, are assumed small enough to be neglected. Between the polarizing and analyzing spin filter, including the sample position, the neutron beam has to be well-collimated throughout the magnetic field in order to maintain polarization. 

In essence $ 3 \times 3$ matrices are measured for each ray, which represent three perpendicular incoming polarization directions with all three directions for analysis. Such measurements can only be recorded by the use of spin turners \cite{3DPNT}. Two flat coils situated before and after the sample turn the spin utilising Larmor precession by $\pi/2$ in two perpendicular directions. By activating none, one or both spin flippers on either side of the sample will allow the measurement of the entire matrix.  

In quantum mechanics, each Pauli matrix is related to an angular momentum operator that corresponds to an observable describing the spin of a spin $\frac{1}{2}$ particle, in each of the three spatial directions.  The three Pauli matrices $\frac{\mri}{2} \tau_1, \frac{\mri}{2} \tau_2$ and $\frac{\mri}{2} \tau_3$ form a basis for the Lie algebra  $\mathfrak{su}(2)$ which exponentiates to $SU(2)$, a double cover for $SO(3)$. Now ${\mri \tau_j}/{2}$ are the generators of a projective representation (spin representation) of the rotation group $SO(3)$ acting on non-relativistic particles with spin $\frac{1}{2},$ such as neutrons, which are fermions. 

An intriguing property of spin $\frac{1}{2}$ particles is that the spin direction must be rotated by an angle of $4\pi$ in order to return to their original configuration. This is due to the two-to-one correspondence between $SU(2)$ and $SO(3)$ and the way that, albeit one imagines spin up or down as the North or South pole on the 2-sphere $S^2$, they are actually represented by orthogonal vectors in $\mathbb{C}^2$. Hence, the states of the particles can be represented as two-component spinors or by $3\times3$ measurement matrices whose entries include the polarization measurement in each direction. 

\def\polzstart{2.5}
\def\polzdepth{2}
\def\polxwidth{1.2}
\def\polyheight{1.5}

\def\polrotzstart{5.25}
\def\polrotzdepth{0.2}
\def\polrotxwidth{1.2}
\def\polrotyheight{1.5}
\def\rotzspacing{0.75}

\def\samplezstart{7.8}
\def\samplelength{1}

\def\anazstart{10.5}
\def\anazdepth{2}
\def\anaxwidth{1.2}
\def\anayheight{1.5}

\def\anarotzstart{9}
\def\anarotzdepth{0.2}
\def\anarotxwidth{1.2}
\def\anarotyheight{1.5}

\def\detzstart{14}
\def\detxwidth{1.2}
\def\detyheight{1.5}
\begin{figure}[!ht]
\begin{center}
\resizebox{0.85\textwidth}{!}{%
\begin{tikzpicture}[y={(-1cm,0.5cm)},x={(1cm,0.5cm)}, z={(0cm,1cm)}]
\coordinate (O) at (0, 0, 0);
\draw[-latex] (O) -- +(1, 0,  0) node [right] {$\xi$};
\draw[-latex] (O) -- +(0,  1, 0) node [left] {$\zeta$};
\draw[-latex] (O) -- +(0,  0, 1) node [above] {$\eta$};

\draw [dashed] (1, 0, 0) -- (\detzstart,0,0);
\draw (0,0,0) node [below] {Neutron Source};

\draw [color=Blue] (\polrotzstart+\rotzspacing,-\polrotxwidth/2,-\polrotyheight/2) -- (\polrotzstart+\rotzspacing,\polrotxwidth/2,-\polrotyheight/2) -- (\polrotzstart+\polrotzdepth+\rotzspacing,\polrotxwidth/2,-\polrotyheight/2) -- (\polrotzstart+\polrotzdepth+\rotzspacing,-\polrotxwidth/2,-\polrotyheight/2) -- (\polrotzstart+\rotzspacing,-\polrotxwidth/2,-\polrotyheight/2) -- cycle;
\draw [color=Blue] (\polrotzstart+\rotzspacing,-\polrotxwidth/2,\polrotyheight/2) -- (\polrotzstart+\rotzspacing,\polrotxwidth/2,\polrotyheight/2) -- (\polrotzstart+\polrotzdepth+\rotzspacing,\polrotxwidth/2,\polrotyheight/2) -- (\polrotzstart+\polrotzdepth+\rotzspacing,-\polrotxwidth/2,\polrotyheight/2) -- (\polrotzstart+\rotzspacing,-\polrotxwidth/2,\polrotyheight/2) -- cycle;
\draw [color=Blue] (\polrotzstart+\rotzspacing,-\polrotxwidth/2,-\polrotyheight/2) -- (\polrotzstart+\rotzspacing,-\polrotxwidth/2,\polrotyheight/2);
\draw [color=Blue] (\polrotzstart+\rotzspacing,\polrotxwidth/2,-\polrotyheight/2) -- (\polrotzstart+\rotzspacing,\polrotxwidth/2,\polrotyheight/2);
\draw [color=Blue] (\polrotzstart+\polrotzdepth+\rotzspacing,-\polrotxwidth/2,-\polrotyheight/2) -- (\polrotzstart+\polrotzdepth+\rotzspacing,-\polrotxwidth/2,\polrotyheight/2);
\draw [color=Blue] (\polrotzstart+\polrotzdepth+\rotzspacing,\polrotxwidth/2,-\polrotyheight/2) -- (\polrotzstart+\polrotzdepth+\rotzspacing,\polrotxwidth/2,\polrotyheight/2);
\draw [color=Blue] (\polrotzstart+\polrotzdepth+\rotzspacing,-\polrotxwidth/2,-\polrotyheight/2) node [right] {  $\pi$/2 spin turner};
\draw [line width = 2pt, color=Blue,-latex] (\polrotzstart+\polrotzdepth/2+\rotzspacing,0,-\polrotyheight/2) -- +(0, 0,  \polrotyheight);

\draw [color=Aquamarine] (\polrotzstart,-\polrotxwidth/2,-\polrotyheight/2) -- (\polrotzstart,\polrotxwidth/2,-\polrotyheight/2) -- (\polrotzstart+\polrotzdepth,\polrotxwidth/2,-\polrotyheight/2) -- (\polrotzstart+\polrotzdepth,-\polrotxwidth/2,-\polrotyheight/2) -- (\polrotzstart,-\polrotxwidth/2,-\polrotyheight/2) -- cycle;
\draw [color=Aquamarine] (\polrotzstart,-\polrotxwidth/2,\polrotyheight/2) -- (\polrotzstart,\polrotxwidth/2,\polrotyheight/2) -- (\polrotzstart+\polrotzdepth,\polrotxwidth/2,\polrotyheight/2) -- (\polrotzstart+\polrotzdepth,-\polrotxwidth/2,\polrotyheight/2) -- (\polrotzstart,-\polrotxwidth/2,\polrotyheight/2) -- cycle;
\draw [color=Aquamarine] (\polrotzstart,-\polrotxwidth/2,-\polrotyheight/2) -- (\polrotzstart,-\polrotxwidth/2,\polrotyheight/2);
\draw [color=Aquamarine] (\polrotzstart,\polrotxwidth/2,-\polrotyheight/2) -- (\polrotzstart,\polrotxwidth/2,\polrotyheight/2);
\draw [color=Aquamarine] (\polrotzstart+\polrotzdepth,-\polrotxwidth/2,-\polrotyheight/2) -- (\polrotzstart+\polrotzdepth,-\polrotxwidth/2,\polrotyheight/2);
\draw [color=Aquamarine] (\polrotzstart+\polrotzdepth,\polrotxwidth/2,-\polrotyheight/2) -- (\polrotzstart+\polrotzdepth,\polrotxwidth/2,\polrotyheight/2);
\draw [color=Aquamarine] (\polrotzstart-\polrotzdepth,-\polrotxwidth/2,-\polrotyheight/2) node [right] { $\;$ $\pi$/2 spin turner};

\draw [line width = 2pt, color=Aquamarine,-latex] (\polrotzstart+\polrotzdepth/2,-\polrotxwidth/2,0) -- +(0,  \polrotxwidth,0);

\draw [color=Sepia] (\polzstart,-\polxwidth/2,-\polyheight/2) -- (\polzstart,\polxwidth/2,-\polyheight/2) -- (\polzstart+\polzdepth,\polxwidth/2,-\polyheight/2) -- (\polzstart+\polzdepth,-\polxwidth/2,-\polyheight/2) -- (\polzstart,-\polxwidth/2,-\polyheight/2) -- cycle;
\draw [color=Sepia] (\polzstart,-\polxwidth/2,\polyheight/2) -- (\polzstart,\polxwidth/2,\polyheight/2) -- (\polzstart+\polzdepth,\polxwidth/2,\polyheight/2) -- (\polzstart+\polzdepth,-\polxwidth/2,\polyheight/2) -- (\polzstart,-\polxwidth/2,\polyheight/2) -- cycle;
\draw [color=Sepia] (\polzstart,-\polxwidth/2,-\polyheight/2) -- (\polzstart,-\polxwidth/2,\polyheight/2);
\draw [color=Sepia] (\polzstart,\polxwidth/2,-\polyheight/2) -- (\polzstart,\polxwidth/2,\polyheight/2);
\draw [color=Sepia] (\polzstart+\polzdepth,-\polxwidth/2,-\polyheight/2) -- (\polzstart+\polzdepth,-\polxwidth/2,\polyheight/2);
\draw [color=Sepia] (\polzstart+\polzdepth,\polxwidth/2,-\polyheight/2) -- (\polzstart+\polzdepth,\polxwidth/2,\polyheight/2);
\draw [color=Sepia] (\polzstart+\polzdepth/2,\polxwidth/2,\polyheight/2) node [left] {Polarizer $\;\;$};
\draw [line width = 2pt, color=Sepia,-latex] (\polzstart+\polzdepth/2,0,-\polyheight/2) -- +(0, 0,  \polyheight);

\draw [color=Plum] (\anazstart,-\anaxwidth/2,-\anayheight/2) -- (\anazstart,\anaxwidth/2,-\anayheight/2) -- (\anazstart+\anazdepth,\anaxwidth/2,-\anayheight/2) -- (\anazstart+\anazdepth,-\anaxwidth/2,-\anayheight/2) -- (\anazstart,-\anaxwidth/2,-\anayheight/2) -- cycle;
\draw [color=Plum] (\anazstart,-\anaxwidth/2,\anayheight/2) -- (\anazstart,\anaxwidth/2,\anayheight/2) -- (\anazstart+\anazdepth,\anaxwidth/2,\anayheight/2) -- (\anazstart+\anazdepth,-\anaxwidth/2,\anayheight/2) -- (\anazstart,-\anaxwidth/2,\anayheight/2) -- cycle;
\draw [color=Plum] (\anazstart,-\anaxwidth/2,-\anayheight/2) -- (\anazstart,-\anaxwidth/2,\anayheight/2);
\draw [color=Plum] (\anazstart,\anaxwidth/2,-\anayheight/2) -- (\anazstart,\anaxwidth/2,\anayheight/2);
\draw [color=Plum] (\anazstart+\anazdepth,-\anaxwidth/2,-\anayheight/2) -- (\anazstart+\anazdepth,-\anaxwidth/2,\anayheight/2);
\draw [color=Plum] (\anazstart+\anazdepth,\anaxwidth/2,-\anayheight/2) -- (\anazstart+\anazdepth,\anaxwidth/2,\anayheight/2);
\draw [color=Plum] (\anazstart+\anazdepth/2,\anaxwidth/2,\anayheight/2) node [left] {Analyzer $\;\;$};
\draw [line width = 2pt, color=Plum,-latex] (\anazstart+\anazdepth/2,0,-\anayheight/2) -- +(0, 0,  \anayheight);

\draw [color=Aquamarine] (\anarotzstart+\rotzspacing,-\anarotxwidth/2,-\anarotyheight/2) -- (\anarotzstart+\rotzspacing,\anarotxwidth/2,-\anarotyheight/2) -- (\anarotzstart+\anarotzdepth+\rotzspacing,\anarotxwidth/2,-\anarotyheight/2) -- (\anarotzstart+\anarotzdepth+\rotzspacing,-\anarotxwidth/2,-\anarotyheight/2) -- (\anarotzstart+\rotzspacing,-\anarotxwidth/2,-\anarotyheight/2) -- cycle;
\draw [color=Aquamarine] (\anarotzstart+\rotzspacing,-\anarotxwidth/2,\anarotyheight/2) -- (\anarotzstart+\rotzspacing,\anarotxwidth/2,\anarotyheight/2) -- (\anarotzstart+\anarotzdepth+\rotzspacing,\anarotxwidth/2,\anarotyheight/2) -- (\anarotzstart+\anarotzdepth+\rotzspacing,-\anarotxwidth/2,\anarotyheight/2) -- (\anarotzstart+\rotzspacing,-\anarotxwidth/2,\anarotyheight/2) -- cycle;
\draw [color=Aquamarine] (\anarotzstart+\rotzspacing,-\anarotxwidth/2,-\anarotyheight/2) -- (\anarotzstart+\rotzspacing,-\anarotxwidth/2,\anarotyheight/2);
\draw [color=Aquamarine] (\anarotzstart+\rotzspacing,\anarotxwidth/2,-\anarotyheight/2) -- (\anarotzstart+\rotzspacing,\anarotxwidth/2,\anarotyheight/2);
\draw [color=Aquamarine] (\anarotzstart+\anarotzdepth+\rotzspacing,-\anarotxwidth/2,-\anarotyheight/2) -- (\anarotzstart+\anarotzdepth+\rotzspacing,-\anarotxwidth/2,\anarotyheight/2);
\draw [color=Aquamarine] (\anarotzstart+\anarotzdepth+\rotzspacing,\anarotxwidth/2,-\anarotyheight/2) -- (\anarotzstart+\anarotzdepth+\rotzspacing,\anarotxwidth/2,\anarotyheight/2);
\draw [color=Aquamarine] (\anarotzstart+\anarotzdepth+\rotzspacing,-\anarotxwidth/2,-\anarotyheight/2) node [right] { $\pi$/2 spin turner};
\draw [line width = 2pt, color=Aquamarine,-latex] (\anarotzstart+\anarotzdepth/2+\rotzspacing,-\anarotxwidth/2,0) -- +(0,  \anarotxwidth,0);

\draw [color=Blue] (\anarotzstart,-\anarotxwidth/2,-\anarotyheight/2) -- (\anarotzstart,\anarotxwidth/2,-\anarotyheight/2) -- (\anarotzstart+\anarotzdepth,\anarotxwidth/2,-\anarotyheight/2) -- (\anarotzstart+\anarotzdepth,-\anarotxwidth/2,-\anarotyheight/2) -- (\anarotzstart,-\anarotxwidth/2,-\anarotyheight/2) -- cycle;
\draw [color=Blue] (\anarotzstart,-\anarotxwidth/2,\anarotyheight/2) -- (\anarotzstart,\anarotxwidth/2,\anarotyheight/2) -- (\anarotzstart+\anarotzdepth,\anarotxwidth/2,\anarotyheight/2) -- (\anarotzstart+\anarotzdepth,-\anarotxwidth/2,\anarotyheight/2) -- (\anarotzstart,-\anarotxwidth/2,\anarotyheight/2) -- cycle;
\draw [color=Blue] (\anarotzstart,-\anarotxwidth/2,-\anarotyheight/2) -- (\anarotzstart,-\anarotxwidth/2,\anarotyheight/2);
\draw [color=Blue] (\anarotzstart,\anarotxwidth/2,-\anarotyheight/2) -- (\anarotzstart,\anarotxwidth/2,\anarotyheight/2);
\draw [color=Blue] (\anarotzstart+\anarotzdepth,-\anarotxwidth/2,-\anarotyheight/2) -- (\anarotzstart+\anarotzdepth,-\anarotxwidth/2,\anarotyheight/2);
\draw [color=Blue] (\anarotzstart+\anarotzdepth,\anarotxwidth/2,-\anarotyheight/2) -- (\anarotzstart+\anarotzdepth,\anarotxwidth/2,\anarotyheight/2);
\draw [color=Blue] (\anarotzstart-\anarotzdepth,-\anarotxwidth/2,-\anarotyheight/2) node [right] { $\;$ $\pi$/2 spin turner};
\draw [line width = 2pt, color=Blue,-latex] (\anarotzstart+\anarotzdepth/2,0,-\anarotyheight/2) -- +(0, 0,  \anarotyheight);

\draw [dashed] (\samplezstart,0,-1.5) -- (\samplezstart,0,1.5);

\draw [color=BrickRed] (\detzstart,-\detxwidth/2,-\detyheight/2) -- (\detzstart,-\detxwidth/2,\detyheight/2) -- (\detzstart,\detxwidth/2,\detyheight/2) -- (\detzstart,\detxwidth/2,-\detyheight/2) -- cycle;
\draw [color=BrickRed] (\detzstart,-\detxwidth/4,-\detyheight/2) -- (\detzstart,-\detxwidth/4,\detyheight/2);
\draw [color=BrickRed] (\detzstart,\detxwidth/4,-\detyheight/2) -- (\detzstart,\detxwidth/4,\detyheight/2);
\draw [color=BrickRed] (\detzstart,-\detxwidth/2,-\detyheight/4) -- (\detzstart,\detxwidth/2,-\detyheight/4);
\draw [color=BrickRed] (\detzstart,-\detxwidth/2,\detyheight/4) -- (\detzstart,\detxwidth/2,\detyheight/4);
\draw [color=BrickRed] (\detzstart,-\detxwidth/2,0) -- (\detzstart,\detxwidth/2,0);
\draw [color=BrickRed] (\detzstart,0,-\detyheight/2) -- (\detzstart,0,\detyheight/2);
\draw [color=BrickRed] (\detzstart,\detxwidth/2,\detyheight/2) node [left] {Detector};

\draw[color=Orange,decorate,decoration={coil,segment length=\samplelength/5*18,amplitude=6}] (\samplezstart,-\samplelength/2,0)   -- (\samplezstart,\samplelength/2,0);
\draw [color=Orange] (\samplezstart,-\samplelength*0.5,0) -- (\samplezstart,-\samplelength*0.5,-1);
\draw [color=Orange] (\samplezstart,\samplelength*0.5,0) -- (\samplezstart,\samplelength*0.5,-1);
\draw [color=Orange] (\samplezstart,0,1)+(-60:.25) [yscale=1.3,->] arc(-60:240:.25);
\draw [color=Orange] (\samplezstart,\samplelength*0.5,0.5) node [left] {Sample};

\end{tikzpicture}
}
\end{center}
\caption{Schematic diagram of measurement apparatus: An ensemble of neutrons are fired from a neutron source which are polarized in the $\eta$-direction by a polarizer. The spin on the neutrons is rotated in the $\zeta$- or $\eta$-directions using a pair of $\frac{\pi}{2}$ spin turners before they reach the sample. Once they traverse the sample another pair of $\frac{\pi}{2}$ spin turners choose the direction of analysis before the analyzer, where neutrons are transmitted along the $\xi$- direction. Finally, the spin measurements are recorded by a position sensitive detector. The sample is rotated around the $\eta$- axis for various tomographic projections. The arrows shown in the polarizer, analyzer and turners are the directions of the magnetic field in those devices. In the case of a spin filter it is the direction of the spin vector that passes through, while in the turners it is the axis about which the spin vector is rotated.}
\end{figure}
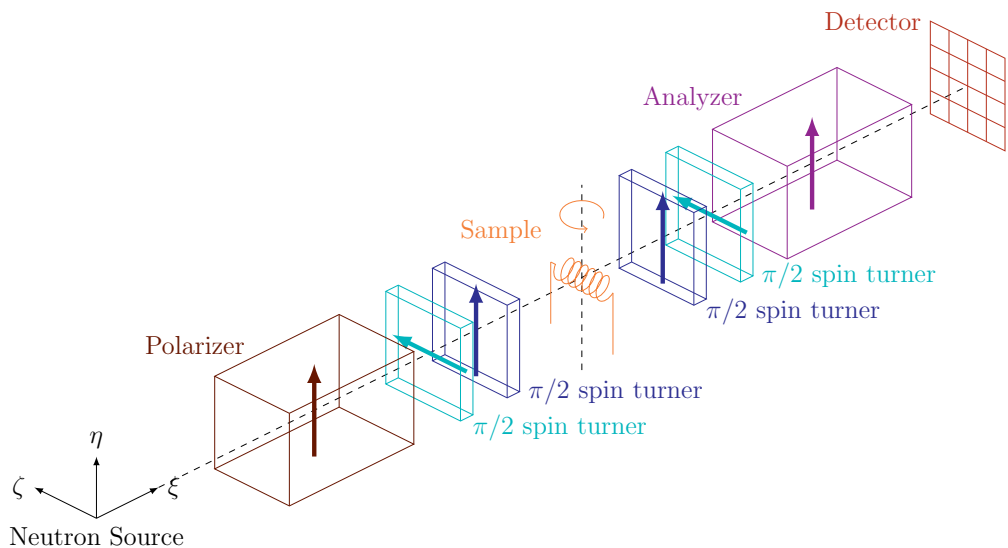

\section{Formulation as a transport equation}

We have formulated the problem in (\ref{eqn:larmor}) as a differential equation for a spin vector along one ray. To apply known results from the mathematical literature we will reformulate this as a transport equation for a matrix valued functions of rays. First we note that the vector product can be written as a matrix vector product
\begin{equation}
\sigma \times B = H(B) \sigma
\end{equation}
where $H(B)$ is a skew symmetric matrix associated with $B$ (the Hodge star)
\begin{equation}
 H(B) =  \left(\begin{array}{rrr} 0 & B_3 & -B_2 \\ -B_3 & 0 & B_1\\ B_2 & -B_1 & 0  \end{array} \right)
\label{eqn:HodgeStar}   
\end{equation}
and so we rewrite (\ref{eqn:larmor}) as 
\begin{equation}
\frac{\mrd \sigma(s)}{\mrd s} = \frac{\gamma_N}{v}H\left( B\left(x\left(s\right)\right)\right)\sigma(s),
\label{eqn:ForwardProblem}
\end{equation}
Along the path of this polarized neutron we have a system of ordinary differential equations (\ref{eqn:ForwardProblem}) which is linear but with variable coefficients. Of course the differential equation cannot be solved simply using the matrix exponential that we would use for a scalar problem as in general the matrix $H(B(x(s)))$ does not commute with its derivative with respect to $s$ (to clarify \cite[eq 3]{TropReview} only holds for uniform fields).

Our physical model assumes that there is no interaction between the direction of travel $\xi$ and the effect of the magnetic field on the neutron. In contrast to many other vector and tensor tomography methods, polarized light tomography for example, the direction of travel $\xi$ of the particle plays no role in its interaction with the quantity being imaged. It plays the same role as the Rytov-Sharafutdinov law in polarized light tomography \cite{SBook}.

The initial condition for this system of ODEs is the specification of $\sigma$ for each ray in the limiting case as $s \rightarrow - \infty$, or as we shall assume $B$ is compactly supported at any point before the ray enters the support of $B$. Considering initial conditions as the three unit basis vectors $e_i$ in turn with $\sigma^i$ as the solution, we assemble a matrix $\Sigma$ whose columns are the spin vectors $\sigma^i$. For any point $x \in \mathbb{R}^3$ and unit vector $\xi$ we have the matrix $\Sigma(x, \xi)$, which now satisfies the generalized transport equation
\begin{equation}
\frac{\partial}{\partial s} \Sigma(x(s), \xi)= \xi \cdot \nabla_x \Sigma(x(s),\xi) =  \frac{\gamma_N}{v} H\left(B(x(s))\right)\Sigma\left(x(s),\xi\right).
\label{eqn:Global_DiffEqn}
\end{equation}
Note that as $H(B)$ is skew symmetric 
\begin{equation}
  \frac{\partial}{\partial s}  \left( \Sigma^T \Sigma \right) =0
\end{equation}
so with the initial condition $ \lim\limits_{s\rightarrow -\infty}\Sigma(x+s \xi ,\xi)=I$ (here $I$ is the $3 \times 3 $ identity matrix) we see that for all $x,\xi \in \mathbb{R}^3, |\xi|=1$ we have  $\Sigma(x, \xi)\in SO(3)$, the group of rotation matrices. Given the transport equation (\ref{eqn:Global_DiffEqn}) we now consider the  `final' data
\begin{equation}
 S(A)(x,\xi) :=  \lim\limits_{s \rightarrow \infty}  \Sigma(x+s\xi,\xi)
\end{equation}
where $A(x) =\frac{\gamma_N}{v} H(B(x))$ as the non-Abelian ray transform data of the matrix valued function $A$. Note that while $S$ is given as a function of $x$ this only serves as an example of the point that the line passes through, it is only the component normal to $\xi$ that is needed.  For rays in one specific plane,
$P_{\eta,z}=\{x\in\mathbb{R}^3: x\cdot \eta=z\}$ normal to a unit vector $\eta$,  
 we will call this $S(A)$  the non-Abelian Radon transform of $A$ restricted to that plane. 

\section{Uniqueness of solution for non-Abelian Radon transforms}

We will now introduce the notation used for a Non-Abelian Radon Transform (NART) in more generality, following \cite{Eskin}, and the references therein.  Rather than the single matrix in  the previous section we now  have three $A_{j}(t)$, $j=0,1,2$ , each a $C^{\infty} \ n\times n$ matrix valued functions on $\mathbb{R}^{2}$, 
with compact support contained in the a ball of some radius $R$. We denote by $\xi = (\xi_1,\xi_2)$  a unit vector and let $\Sigma(x,\xi)$ be the matrix solution of the partial differential equation
\begin{equation}
\xi \cdot \nabla_x{ \Sigma(x,\xi)} = (A_1(x)\xi_1 + A_2(x)\xi_2 + A_0(x))\Sigma(x,\xi),
\label{eqn:General_eqn}
\end{equation}
with the initial condition $\Sigma(x + s\xi,\xi) \rightarrow I_n,$ as $s \rightarrow -\infty$. 
We then define the limit of $\Sigma(x + s\xi,\xi) $ as $s \rightarrow \infty$ to be the NART of $A(x,\xi)=A_1(x)\xi_1 + A_2(x)\xi_2 + A_0(x)$, which we denote by $S(A)$.

 Note we have allowed a first order dependence of the right hand side on the direction $\xi$ which appears in other physical and geometric problems.
 This case does not include polarized light tomography  in which there is a 4th order dependence on $\xi$. This is apparent from the Rytov-Sharafutdinov law defined in \cite[\textsection 5.1.4]{SBook}.

\begin{definition}
The matrices $A^{(1)}_j(x)$, $0\le j \le 2$ and $A^{(2)}_j(x)$,$0 \le j \le 2$ are said to be {\em gauge equivalent} if 
there exists a nonsingular $n\times n$ $C^\infty$ matrix valued function $g(x)$ such that $g(x)=I_n$ for $|x|\ge R$ and 
\begin{eqnarray}
A^{(2)}_j &=& gA^{(1)}_j g^{-1} + \frac{\partial{g}}{\partial x_j}  g^{-1},	j = 1, 2,\\
A^{(2)}_0 &=&gA^{(1)}_0g^{-1}.
\end{eqnarray}
\label{def:GaugeEquivalent} 
\end{definition}
The significance of this is that gauge equivalent $A^{(1)}(x,\xi)$ and $A^{(2)}(x,\xi)$ have the same NART, 
so we can can hope at best that the inverse problem has a unique solution up to gauge equivalence.
Indeed \cite{Eskin} proves the following theorem which depends on holomorphic matrix integrating factors. 
\begin{theorem}(Eskin)\label{thm:gaugeunique}
Suppose	$A^{(1)}_j(x)$ and $A^{(2)}_j(x)$, $0\le j \le 2$, are $C^\infty$ compactly supported matrices with 
the same NART. Then $A^{(1)}_j(x)$, $0\le j \le 2$ and $A^{(2)}_j(x)$, $0\le j \le 2$ are gauge equivalent.
\end{theorem}
There is also a more general geometric formulation of this problem on a manifold where $A_1$ and $A_2$ are components of connection and $A_0$ is regarded as a Higgs field. For dimension 2 \cite{paternain2012attenuated} proves the equivalent of Theorem \ref{thm:gaugeunique} in this context. The proof is somewhat simpler in that it uses only scalar holomorphic integrating factors. 

For less smooth $A$ there are currently only local uniqueness results, that is the data determines $A$ if it is known a priori that some norm of $A$ is small enough. Novikov \cite[Corr. 5.3]{Novikov} proves such a result using the norm $||A||_{\alpha,1+\epsilon,\rho}$ defined by \cite[eq 2.4]{Novikov} as Holder continuous of order $\alpha>0$ while as $|x|\rightarrow \infty$ and the decay of $A$ and its Holder quotient is at least $O(|x|^{-1-\epsilon})$, $\epsilon>0$. Novikov's result shows that provided $||A||_{\alpha,1+\epsilon,\rho}<C(\alpha,\epsilon,\rho)$ for some constant $C$ the solution to the inverse problem is unique up to a gauge transformation.

\section{Uniqueness of solution for $B$}
\label{section:ForwardProblemFormulation}

We will apply Theorem \ref{thm:gaugeunique} to the case of recovery of $B$ restricted to some plane from the polarimetric neutron data for all rays in that plane $P_{\eta,z}$. We assume that the restriction of $B$ to $P_{\eta,z}$ has compact support. For notation  convenience we identify $P_{\eta,z}$ with $\mathbb{R}^2$ and our problem on the plane,  in the notation of the non-Abelian Radon transform, has $n=3$, $A_1=A_2=0$ and $A_0=\frac{\gamma_N}{v}H(B(x))$.

Gauge equivalence given $A_1=A_2=0$ implies that $\partial g/\partial x_j=0 $ for $j=1,2,$ which means that $g$ must be constant and hence the identity. This means that for $A(x,\xi) = \frac{\gamma_N}{v} H(B(x(t)))$ gauge equivalence means equality.
 Now Theorem~\ref{thm:gaugeunique} yields the following
\begin{theorem}Given a vector field $B$ such that $B|_{P_{\eta,z}}$ is smooth with compact support, the data $S((\gamma_N/v)H(B))(x,\xi)$ for $x\in P_{\eta,z}$, $\eta\cdot\xi=0$ uniquely determines  $B|_{P_{\eta,z}}$.
\label{thm:PlaneBReconFormula}
\end{theorem}

There are two important limitations here. Firstly we assumed that $B$ is compactly supported on each plane where we  measure. As there is no perfect magnetic shielding of course this is not physically accurate. 

Secondly, while on a sufficiently small scale $B$ will be smooth, in practice we will represent the field as samples on a grid, and it may be that there are large jumps in the magnetic field on that scale. This means that if Theorem \ref{thm:PlaneBReconFormula} breaks down and for less smooth fields the solution is non-unique that may create practical problems for solving the inverse problem.

For a sufficiently large $B$  or small $v$ what is informally known as `phase wrapping' might occur along a ray path \cite{TropReview}. For example if Euler angles are used as a coordinate chart on $SO(3)$ one of the coordinates might change by more than $\pi$ along the ray. As pointed out in \cite{TropReview}, for data from one direction $\xi$ one cannot know if the spin vector has completed several turns around some axis  along the path. Theorem \ref{thm:PlaneBReconFormula} shows that when all rays are measured through the support of $B|_{P_{\eta,z}}$ the magnetic field, and indeed the spin along the rays, is uniquely determined on that plane. So this phase wrapping does not result in any ambiguity. As we shall see it does complicate the reconstruction. On the other hand an advantage of using a polyenergetic neutron source and time of flight measurements is that data is a available for a range of values of $v$, and the largest value available might result in no phase wrapping for a given  $B$.

For Holder continuous $B$ that are not smooth, Novikov's theorem guarantees uniqueness of solution $B$ provided the norm of  $B/v$ (on a given plane) is small enough. For fixed $B$ there will be a $v$ such that uniqueness is guaranteed by this result, so we can recover $B$ on this plane provided we make a measurements with sufficiently energetic neutrons.

 \section{Linearization}

For convenience let us now define the forward map as 
\begin{equation}
 \mathcal{S}(B) =  S(\frac{\gamma_N}{v}H(B))
\end{equation}
and notice that this is non-linear. We seek a linearization, that is a Gateaux derivative, of this map about a fixed $B$.


Note that the transpose of $\Sigma$ satisfies
\[
\frac{\partial}{\partial s} \Sigma^T = - \frac{\gamma_N}{v}H(B)\Sigma^T
\]
where we have suppressed the arguments $(x+s \xi,\xi)$ for brevity.
Now consider a perturbation in which $B$ replaced by $B+\delta B$ and $\Sigma$ by $\Sigma+\delta \Sigma= \mathcal{S}(B+\delta B)$. Of course $\lim\limits_{s\rightarrow -\infty}\delta\Sigma=0$.
Substituting in (\ref{eqn:Global_DiffEqn}) and ignoring second order terms
\[ \frac{\partial}{\partial s} \delta\Sigma =\frac{\gamma_N}{v}\left( H(B) \delta \Sigma  + H(\delta B)  \Sigma \right) +O(\delta B)^2.
\]
Now consider
\[ 
 \frac{\partial}{\partial s}\left( \Sigma^T \delta \Sigma \right) =   \Sigma^T   \frac{\partial}{\partial s}\delta\Sigma  +  \frac{\partial}{\partial s} \left( \Sigma^T \right) \delta \Sigma +O(\delta B)^2
\]
using the expressions above for the derivatives of $\delta \Sigma$ and for $\Sigma^T$
\begin{equation}\label{equ:ATRT}
\frac{\partial}{\partial s}\left( \Sigma^T \delta \Sigma \right) =\frac{\gamma_N}{v} \Sigma^T H(\delta B) \Sigma +O(\delta B)^2
\end{equation}
so that
\[
\lim\limits_{s\rightarrow\infty}\left(\Sigma^T \delta \Sigma \right) =\frac{\gamma_N}{v} \int\limits_{-\infty}^\infty \Sigma^T H(\delta B) \Sigma \, \mrd s +O(\delta B)^2.
\]
Now we notice that $\Sigma^TH(B)\Sigma = H(\Sigma^T B)$ as applying a rotation is just a change of coordinates and this is the transformation rule for a tensor and vector respectively. So we see that  a simpler  formula for the perturbed data is
\[
\left(\Sigma^T \delta \Sigma \right)_{s\rightarrow\infty} = \frac{\gamma_N}{v}H\left( \int\limits_{-\infty}^\infty \Sigma^T\delta B\, \mrd s\right) +O(\delta B)^2.
\]
Going back to (\ref{equ:ATRT}) we see this is an attenuated ray transform in the terminology of \cite{paternain2012attenuated}, with a skew-symmetric Higgs field and flat connection. From that paper we see that a smooth compactly supported $\delta B$ on a plane is determined uniquely by  $D\mathcal{S}_B \left( \delta B\right)$ on that plane, where  the derivative D$\mathcal{S}$ is defined by
\begin{equation}
D\mathcal{S}_B \left( \delta B\right) = \frac{\gamma_N}{v}\mathcal{S}(B) H\left( \int\limits_{-\infty}^\infty \Sigma^T\delta B\, \mrd s\right).    
\end{equation}
In other words the derivative is injective.

For a field $B$ that is small in the plane the (so in the above replacing $B$ by zero and $\delta B$ by $B$) linear approximation gives us
\begin{equation}
e_1\cdot \mathcal{S}(B)(x,\xi) e_2 = \frac{\gamma_N}{v}X(B_3)(x,\xi) +O(|B|)^2.
\label{eqn:3Data}
\end{equation}
where $X$ is the scalar Radon transform in the plane.
With cyclic permutations of the indices $(1,2,3)$ it is clear we can retrieve the magnetic field. This can be done using any two-dimensional Radon inversion method, and this agrees with what is already seen in the experimental literature  e.g.  \cite{Kardjilov,3DPNT}.

\section{Forward problem simulation}
\label{section:forward_simulation}

\begin{figure}[!ht]
\centering
\includegraphics[scale = 0.37]{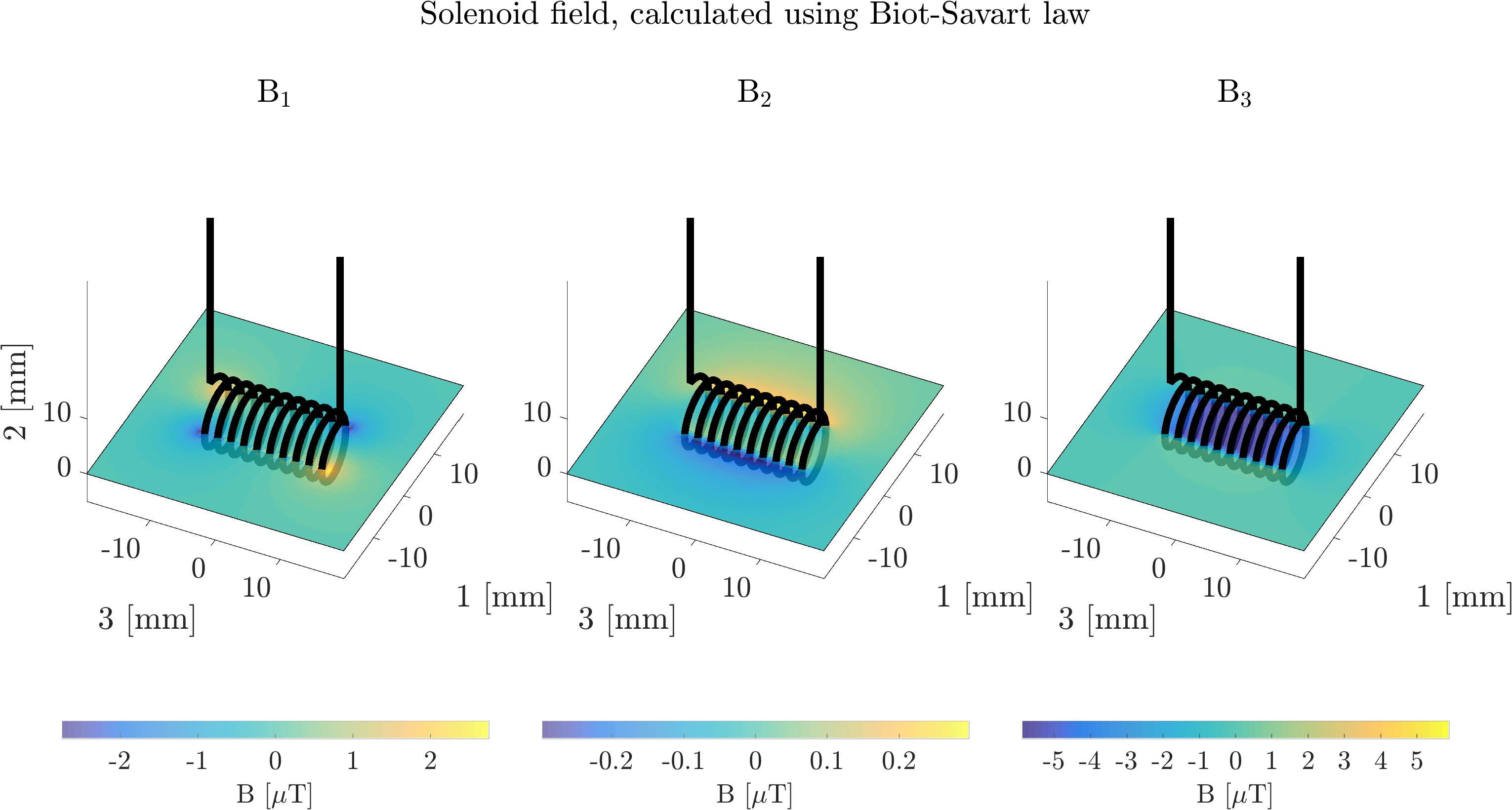}
\caption{Central slice of a simulated solenoid showing components of the magnetic field}
\label{fig:central_solenoid}
\end{figure}

Both to create simulated data ans as part of a solution method we need to solve the forward problem numerically for known magnetic fields.
To create an interesting test field, and one that is practical for experimental tests we assume a constant current passes through and one-dimensional wire and calculate the magnetic field using the Biot-Savart law, performing numerical integration along a piecewise straight continuous curve, The central slice of a simulated solenoid is shown in Figure \ref{fig:central_solenoid}. 
The magnetic field  is calculated at the central point of each voxel in a uniform grid.

Let us look at the procedure to generate simulated PNTMF data. Detail of how the forward solver operates can be explained by taking a single ray at a time. Consider a ray emanating from a point $x_0$ outside the voxel grid, in direction $\xi$, so that $\xi$ is normal to the detector plane and the ray passes through the centre of a detector pixel.  We set  $\Sigma(x_0,\xi) = I$. We assume that the field $B$ is uniform on each voxel and zero outside the grid. We calculate the entry and exit point of the ray along the voxel grid. In conventional tomography one typically implements a subroutine that returns the length of intersection of a ray with a cubic voxel, and this forms a sparse matrix used in forward and inverse calculation. We used the method of Jacobs \cite{jacobs} as implemented in \cite{DS} for tensor tomography problems for this calculation, with the modification that the order in which the voxels are encountered by the ray is also returned by the subroutine. We then solve the ODE (\ref{eqn:ForwardProblem}) along this ray, noting that on voxels where $B$ is constant we have an analytical solution

Consider the eigenvalues of the skew symmetric matrix $(\gamma_N/v)H(B(x))$ for fixed $x$ and  $B(x)\ne 0$ 
\begin{equation}
\lambda_1 = 0, \ \lambda_2 = \mri(\gamma_N/v)|B|, \ \lambda_3 = -\mri(\gamma_N/v)|B|,
\label{eqn:EigsM}
\end{equation}
where the the null space is spanned by $B$.
The solution to (\ref{eqn:ForwardProblem}) for part of a ray on which $B$ is constant takes the
\begin{equation}
\sigma(s) = c_1 B + c_2 \cos (\gamma_N|B|s/v) w_+ + c_3 \sin (\gamma_N|B|s/v)w_-,
\label{eqn:SolnFormulation1}
\end{equation}
for real vectors $w_\pm$  orthogonal to $B$ and constants $ c_1,c_2,c_3 \in \mathbb{R}$ determined by initial conditions on entering the voxel. Note that the component in the same directions as $B$ remains fixed while the components orthogonal rotate.

Rodrigues' rotation formula gives a simple expression for the rotation matrix about a unit vector $k$ anticlockwise through an angle $\phi$ 
\begin{equation}
 \mathrm{R}(k,\phi) = I + (\sin\phi)H(k) + (1-\cos\phi)H^2(k),
 \label{eqn:MatrixRodrigues}
\end{equation}
For a given ray through $x_0$ let $x_i = x_0 + s_i \xi$ be the point at which the ray enters the $i$th voxel on its path and ${B}^i$ the constant value of $B$ assumed on that voxel. Then
\[
\Sigma(x_0+ s_{i+1}\xi,\xi) = \mathrm{R}\left(\frac{B}{|B|}, \frac{\gamma_N|B|(s_{i+1}-s_i)}{v}\right) \Sigma(x_0+ s_{i}\xi,\xi)
\]
and hence our numerical approximation to $S(B)(x_0,\xi)$ is
\[
\prod\limits_{i} \mathrm{R}\left(\frac{B}{|B|}, \frac{\gamma_N|B|(s_{i+1}-s_i)}{v}\right)
\]
and this is repeated for each ray path, which in our experiment means that we have an $x_0$ for each detector pixel and $\xi$ is rotated at equiangular increments. 
The simulated data have been validated experimentally, as shown in of \cite[Figure 2]{3DPNT} where experimental PNTMF data matches our simulated data.

\section{Reconstructing the magnetic field from PNTMF Data}
\label{section:ReconPNMFTLinear}

The central slice of the solenoid (magnetic field) simulated by the use of the Biot-Savart law in Figure \ref{fig:central_solenoid}, is utilized in the forward model by the process described above to generate data. Initially the solenoid was simulated on a $180\times180$ pixel grid. 
To simulate the data 270 rays of neutrons (uniform velocity with wavelength of $ 5$~\AA, speed 790 m/s) were fired for every angular increment (1 degree in this case) of the usual parallel beam tomographic data acquisition process. The data is rebinned by a factor of three to give the data which is three sets of $90 \times 120$ arrays for each of the three components from the spin matrix. Furthermore $5\%$ Gaussian pseudo-random noise was added.

\begin{figure}[!ht]
\centering
\includegraphics[scale=0.75]{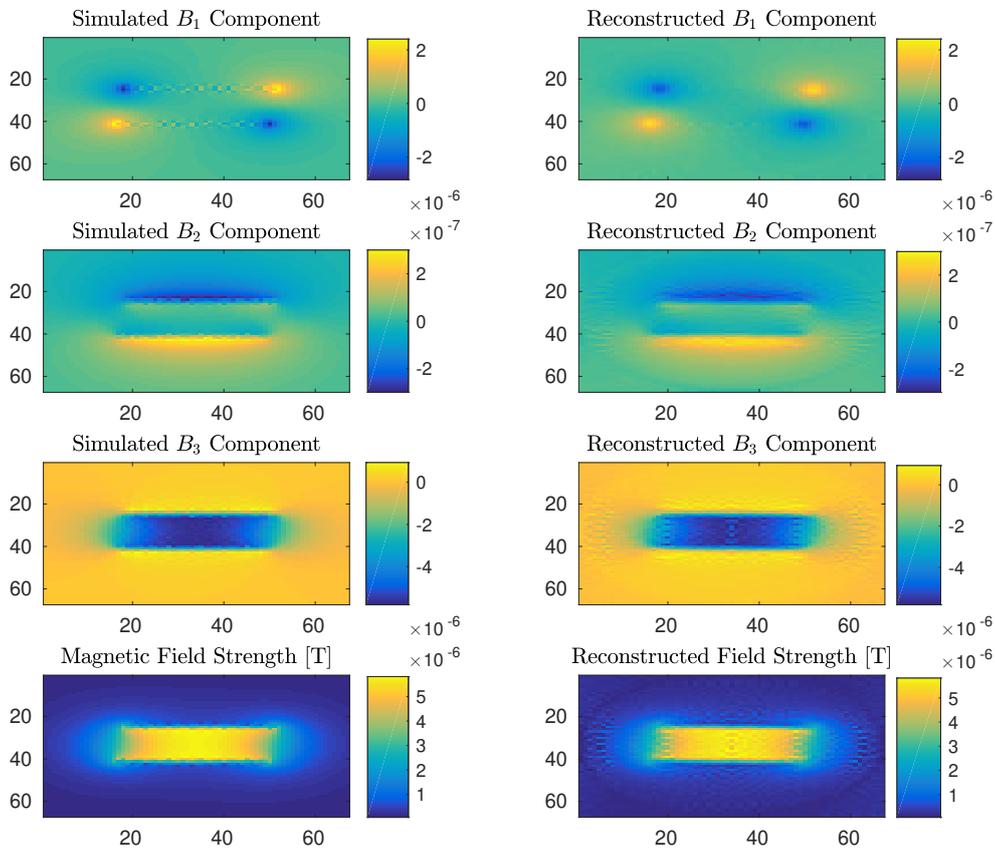}
\caption{The reconstruction of PNTMF data of a central slice of a weak solenoid.}
\label{fig:LinearPaperResults}
\end{figure}

First we implement the simple reconstruction using the linear approximation about $B=0$.
Radon transform inversion was implemented using first a discrete Fourier transform to implement a  a Hamming filter. The filtered data is then backprojected onto a $67 \times 67$ pixel grid to achieve the reconstructed components of the magnetic field. The relative errors from top to bottom in Figure \ref{fig:LinearPaperResults} corresponding to $B_1$, $B_2$, $B_3$ and $|B(x(t))|$ (magnetic field strength) are $20\%, 16\%, 11\%$ and $9\%$ respectively. Since reconstruction is of the central slice of a solenoid, with strength approximately $5.8 ~\mu T,$ the maximum a single neutron precesses as it passes through the domain is 2 degrees. This is well within the range for the linearized problem to work.  

\begin{figure}[!ht]
\centering
\includegraphics[scale=0.7]{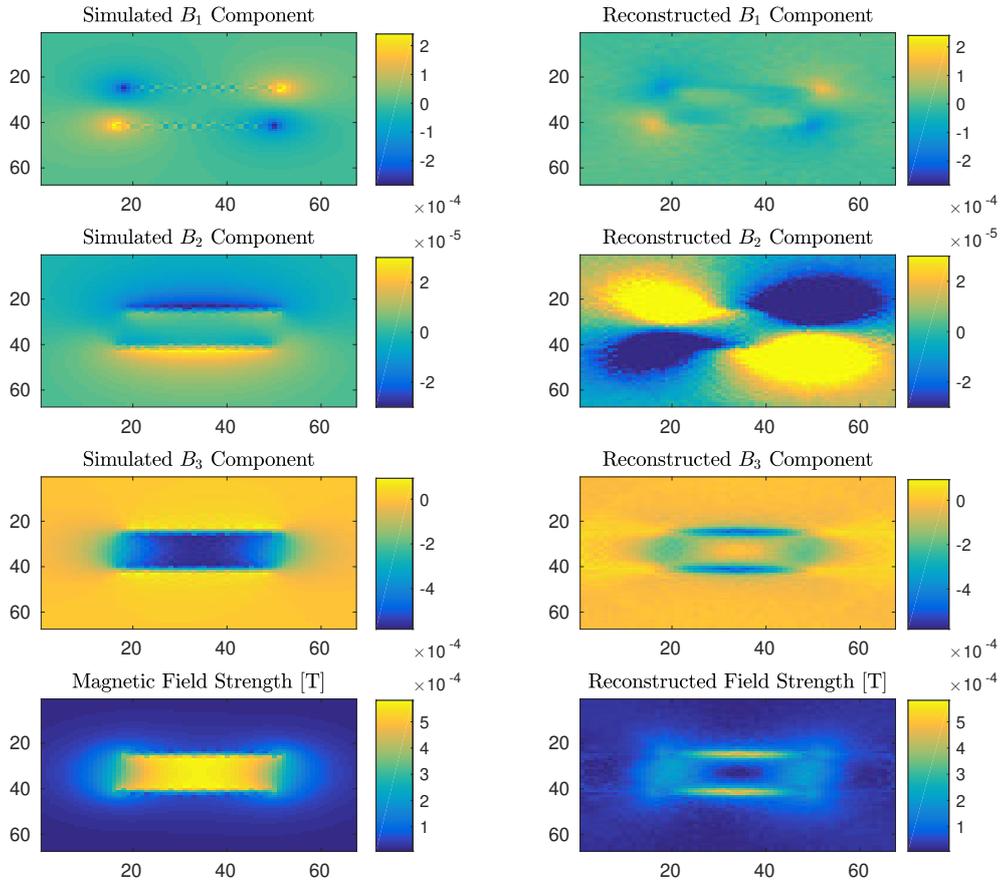}
\caption{Method of (\ref{eqn:3Data}) fails for stronger magnetic fields.}
\label{fig:LinearPNMFTFail}
\end{figure}

When the strength of the magnetic field increases to the extent that a single neutron precesses more than $14^{\circ}$ approximately, this specific method fails. This is when the small angle approximation breaks down, i.e. when $\sin\phi$ differs significantly from  $\phi.$ One such illustration is present in Figure \ref{fig:LinearPNMFTFail} where one notices the artifacts coming through. The relative error for the $B_2$ component is $1.5$

\section{Nonlinearity}
\label{section:NonlinearIP}

\begin{figure}[!ht]
\centering
\includegraphics[scale=0.9]{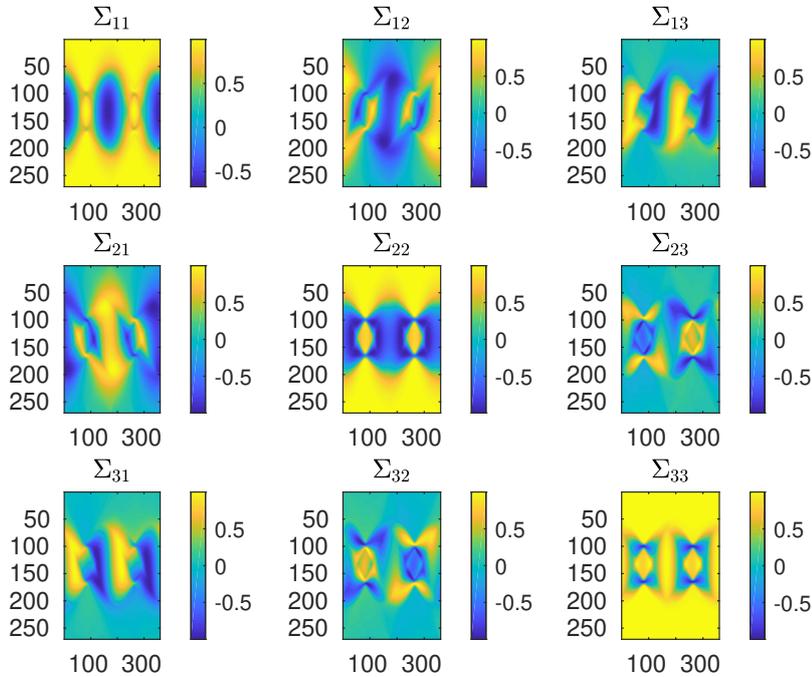}
\caption{Spin data displayed as a matrix of sinograms for $\Sigma$  for $\mathcal{S}(200 B)$ where $B$ is the solenoid field used before.  Contrast with Figure \ref{fig:NonlinearSepTest} to see non-linearity in data.}
\label{fig:NonlinearCombTest}
\end{figure}

\begin{figure}[!ht]
\centering
\includegraphics[scale=0.9]{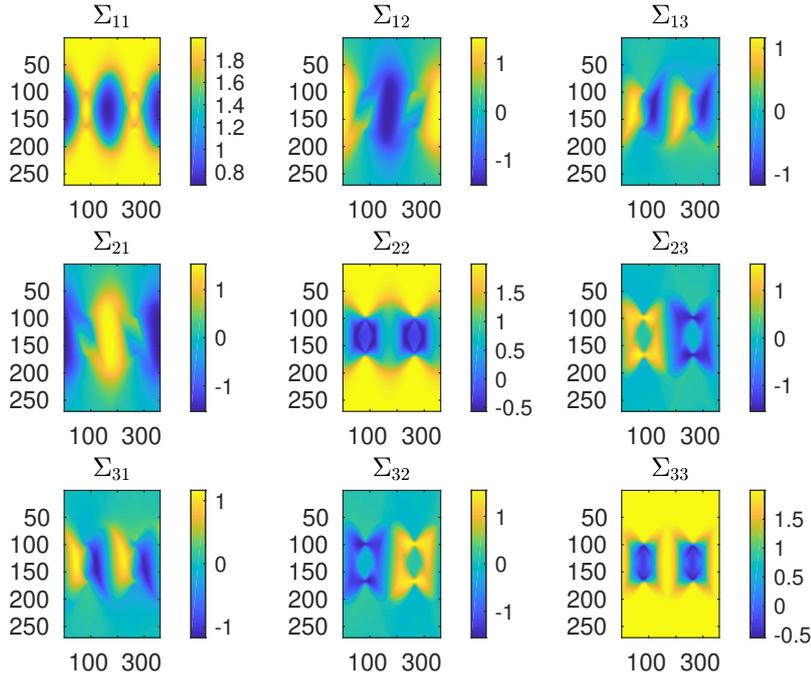}
\caption{By contrast to Figure \ref{fig:NonlinearCombTest} this is a plot of  $\mathcal{S}(50 B) + \mathcal{S}(150 B)$. Note this differs from $\mathcal{S}(50 B+150 B)$}
\label{fig:NonlinearSepTest}
\end{figure}

Non-linear inverse problems have a more complex connection between data and model. For PNTMF, this means   $\mathcal{S}$ is  non-linear operator. When testing a problem for non-linearity practically it is useful to consider the failure of scaling property $\mathcal{S}(\alpha B) = \alpha \mathcal{S}(B)$ and the superposition property 
$\mathcal{S}( B^1 + B^2) = \mathcal{S}(B^1) +  \mathcal{S}(B^2)$ separately. In Figures \ref{fig:NonlinearCombTest}  and \ref{fig:NonlinearSepTest},  we compare the values of each of the 
nine sinograms in the data set obtained scaling the magnetic field is the central slice of the solenoid utilized. 
This linear approximation is seen to break down when the spin vector rotates significantly about some axis. The spin on the neutron depends upon two things; the strength of magnetic field and the time spent by the neutron in the magnetic field. The linear approximation fails and therefore we have to resort to methods in which we are able to solve non-linear problems in order to successfully image magnetic structures in magnetic materials.  

We can solve the inverse problem as an optimization problem: we seek to  minimize $||\mathcal{S}(B)- \mathcal{S}_{\mathrm{meas}} ||^2 $, where $\mathcal{S}_{\mathrm{meas}}$ is the measured data. As $\mathcal{S}(B)$ is non-linear it is possible that the function being minimized is not-convex. We can see this non-convexity in just one direction $B$ by plotting $||\mathcal{S}(\alpha B)- \mathcal{S}_{B} ||^2$ as the parameter $\alpha$ varies and noting if the curve is convex. We see in Figure \ref{fig:non-convex}  convexity breaks down for $\alpha$ over a large enough range.

\section{Modified Newton Kantarovich Method}
\label{section:MNKM}

In general Newton-Kantarovich methods solve a non linear problem by successively solving a linearization of the problem and updating the solution. In the context of inverse problems it is common to apply these methods to overdetermined problems and to employ some regularization in the inversion. Gradient decent methods are related optimization methods that seek a critical point of an objective function such as the squared error between the simulated and measured data. The repeated solution of the linearized problem is often costly computationally so a variation is to to use the linearaization about some fixed value (in our case zero) for which we have an explicit inverse for the linearized problem. When the derivative at a fixed value is used, following \cite[Ch 2]{NewtonK} we call such methods Modified Newton Kantarovich method(MNKM).
Typically the convergence of such methods requires more iteration and the radius of convergence smaller compared to the full method.  Suppose the true magnetic field is $B_{\mathrm{true}}$ and the data $\mathcal{S}_{\mathrm{meas}}= \mathcal{S}(B_{\mathrm{true}}) + \mathrm{errors}$. So we seek a solution $B$ that results in a small residual $||\mathcal{S}(B)- \mathcal{S}_{\mathrm{meas}} ||^2 $. To this end we start with an initial approximation, typically $B^0=0$ (superscripts are now iteration numbers) and at each iteration to derive an update 
\begin{equation}
B^{n+1} = B^{n} + \delta B^{n}, \quad n = 0,1,2,....,
\label{eqn:NKScheme}
\end{equation}
where $B$ is the magnetic field desired and $\delta B$ is the update satisfying the linear operator equation
\begin{equation}
D\mathcal{S}_{B^{0} }\delta B^{n} = \mathcal{S}_{\mathrm{meas}}-\mathcal{S}(B^{n}).
\label{eqn:JacobianEveryStep}
\end{equation}

A damped MNKM method has a line search parameter, $\alpha,$ which controls the extent of the update, given as  
\begin{equation}
B^{n+1} = B^{n} + \alpha ~ \delta B^{n}, \quad n = 0,1,2,....
\label{eqn:MNKM}
\end{equation}
A line search is performed at every step to to choose the a good  $\alpha$ for the update direction $\delta B^{n}$. The forward problem s solved three times, each time with a different parameter, $\alpha$. A quadratic is fitted to the residual error as a function of $\alpha$ and the $\alpha$ value for which the residual is at the minimum, is chosen as the update parameter.

\begin{algorithm}[H]\label{algorithm:MNKMethod}
\SetAlgoLined
\KwData{For a given plane the PNTMF data $\mathcal{S}$ for a fixed $v$ and for a set of parallel rays at equiangular increments, maximum number of iterations $MaxIt$  and convergence tolerance $\mathrm{TOL}$}
\KwResult{Approximation to three components of  $B$ on a voxel grid on that plane}
 $B^0 = 0$\;
 $n = 1$\;
\Repeat{ $(B^{n} - B^{n-1}) < \mathrm{TOL}$ or $n>MaxIt$}{
$\delta B_{i}^{n} = X^{-1}\left[ e_j \cdot [\Sigma - \mathcal{S}(H(B^n))] \cdot e_k \right]$\;
Line search to find $\alpha$ \;
$B^{n+1} = B^{n} + \alpha(\delta B^{n})$\; 
$n = n+1$\;
}
\caption{Modified Newton Kantarovich Method}
\end{algorithm}

In order to test the MNKM the solenoid in Figure \ref{fig:central_solenoid} is scaled up by a factor of 50 which means the magnetic field is of strength $290 ~\mu$T. Using the implementation of the forward model in Section \ref{section:forward_simulation}, initial data is generated for the stronger magnetic field. 
Thereafter the reconstruction process adopted in the linearized inversion process is utilized to yield the result from the first iterate. This is fed back into the forward model to obtain new data which is subtracted from the initial data to give the data set on which the reconstruction procedure will be employed. Upon completion of Radon inversion, the update, $\delta B$, is found.. The iterative procedure terminates at a predetermined tolerance, TOL $= 10^{-5}$.

\begin{figure}[!ht]
\centering
\includegraphics[scale=0.5]{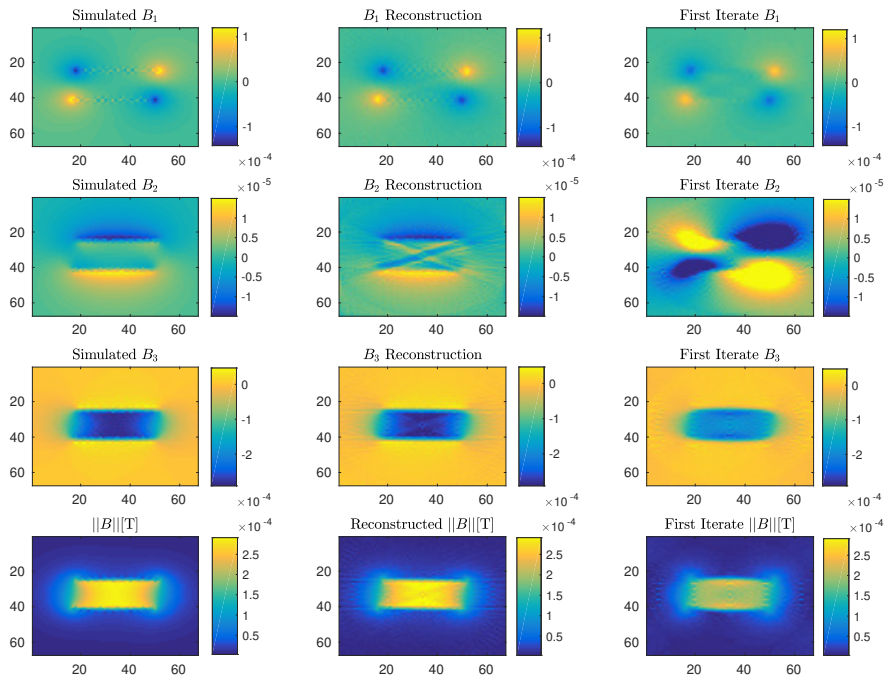}
\caption{Convergence of iterative procedure.  Columns show actual simulated field, reconstruction and first iteration. Rows are the three components of the $B$ field and the norm.}
\label{fig:MNKMConverges}
\end{figure}

Figure \ref{fig:MNKMConverges} shows the results for the iterative reconstruction  Algorithm \ref{algorithm:MNKMethod}. The initial spin data was simulated using a $180\times180$ pixel grid where 
270 rays of neutrons (uniform velocity with wavelength $ = 5$~\AA) were fired for every angular increment (1 degree in this case) of the usual parallel beam tomographic data acquisition process. The data is binned by a factor of three to give the data which is three sets of 
$90 \times 120$ arrays. Furthermore, $5\%$ pseudo-random Gaussian noise was added. Reconstruction was performed on a grid which 
does not evenly divide the grid used for simulation, i.e. $67 \times 67$. The relative errors are $22\%, 17\%, 9\%$ and $10\%$ 
for the magnetic field strength, $\vert B(x(t)) \vert$ and components $B_1, B_2$ and $B_3$ respectively. 
It took 25 iterates to converge and the maximum a neutron precesses throughout this specific magnetic field is $88^{\circ}$. The amount of precession can be calculated using the  forward solver.

\begin{figure}[!ht]
\centering
\includegraphics[scale=0.5]{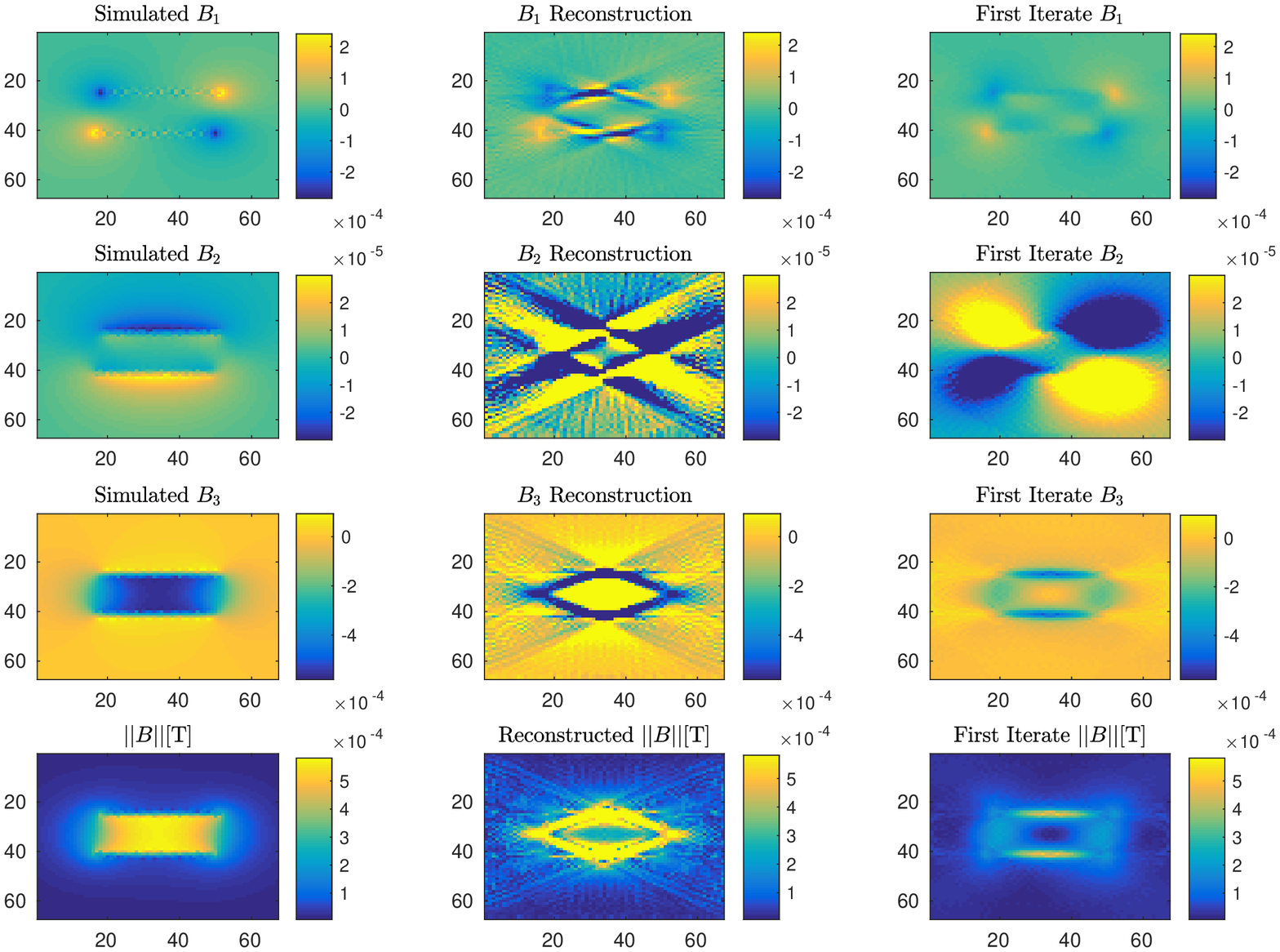}
\caption{Compared to Figure \ref{fig:MNKMConverges}, the method fails to converge with a  magnetic field $100$ times larger.}
\label{fig:MNKMDiverges}
\end{figure}

However, one drawback of this method is illustrated when the solenoid in Figure \ref{fig:central_solenoid} is scaled up by a factor of 100 which means the magnetic field is of strength of $580 ~\mu$T. In this situation the line search in the  MNKM gets trapped and makes no progress reducing the residual.  
Figure \ref{fig:MNKMDiverges} illustrates this where the magnetic field is increased by a factor of $100$.

\begin{figure}[!ht]
\centering


\includegraphics[width=0.47\textwidth]{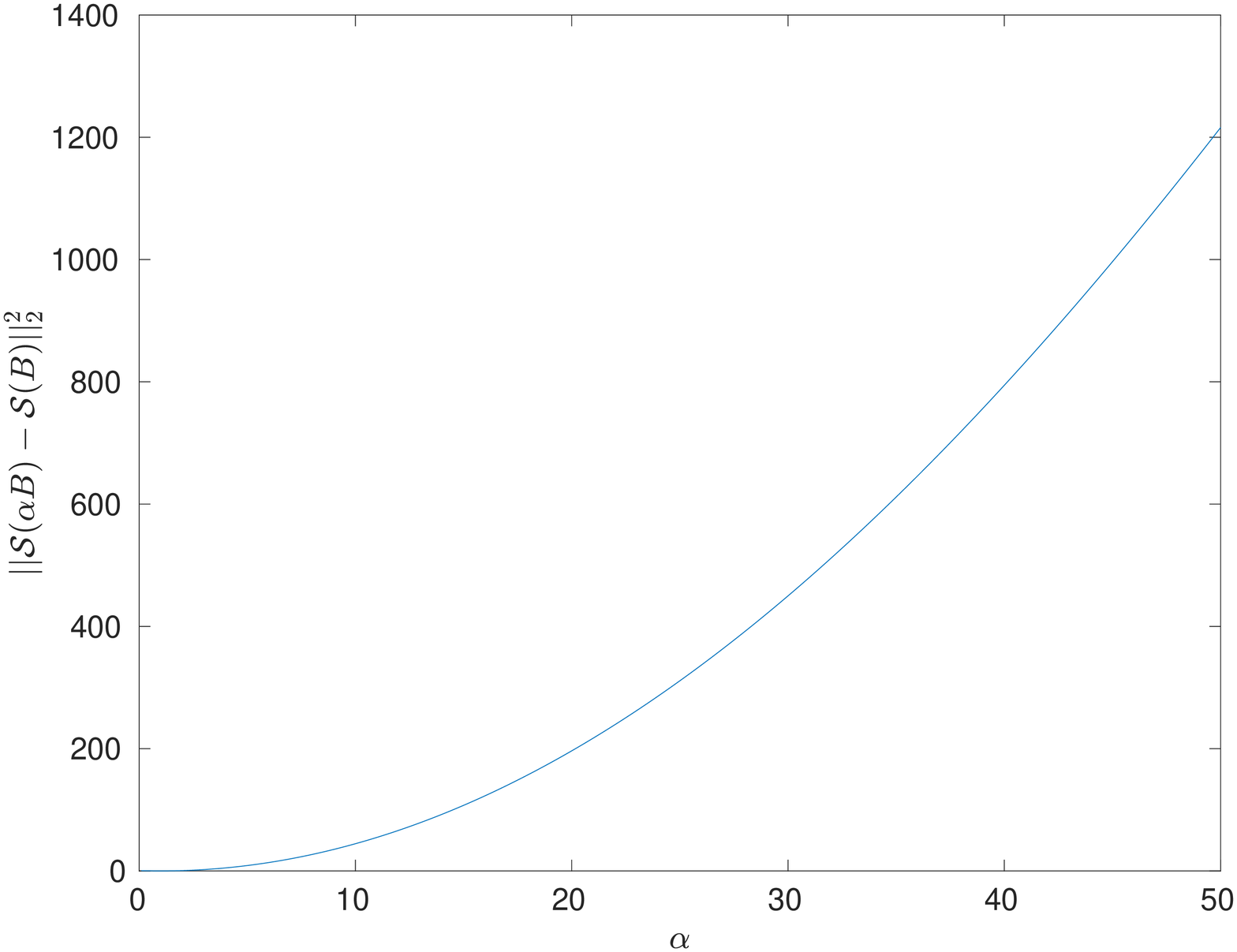}
\includegraphics[width=0.47\textwidth]{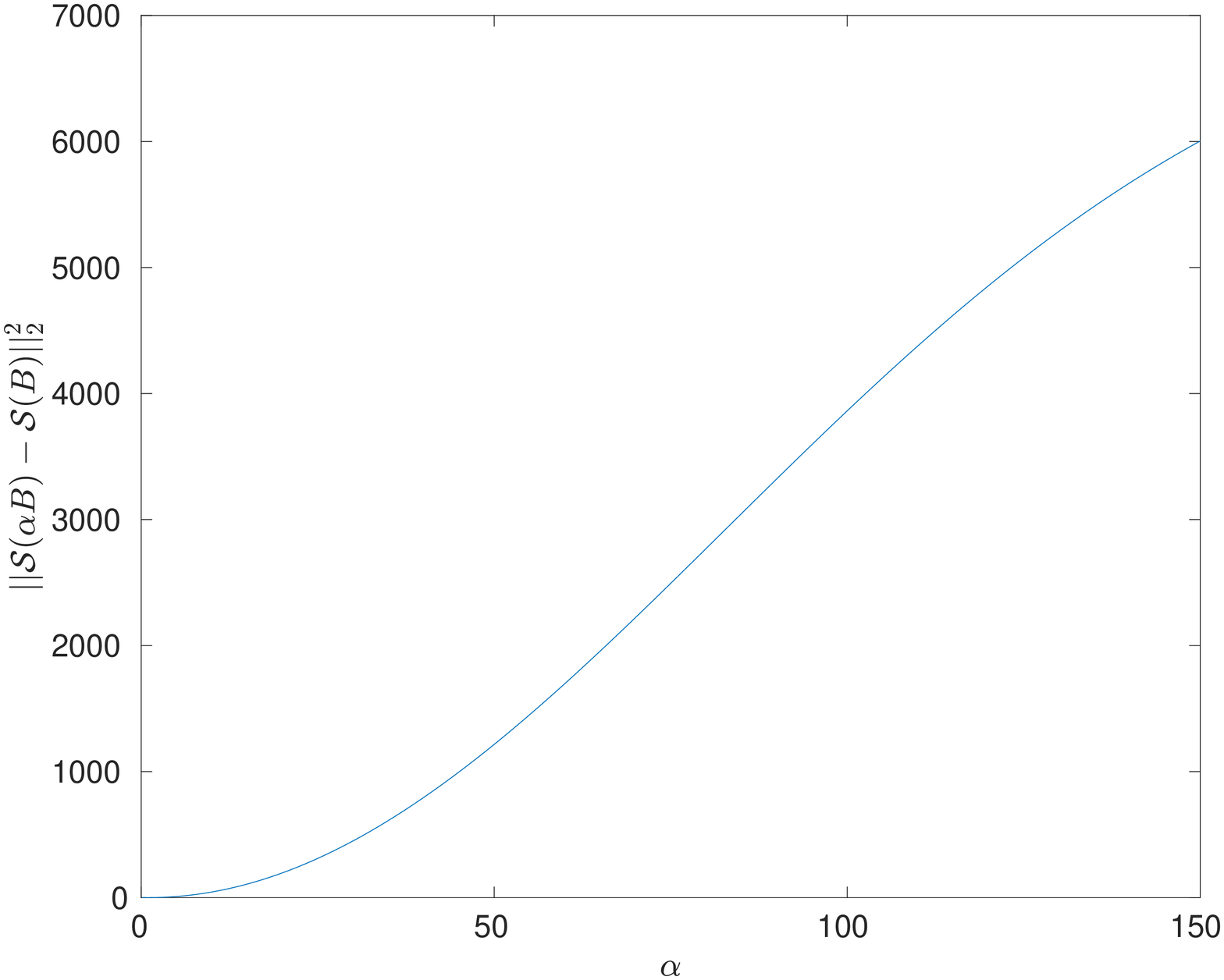}

\centering
\includegraphics[width=0.47\textwidth]{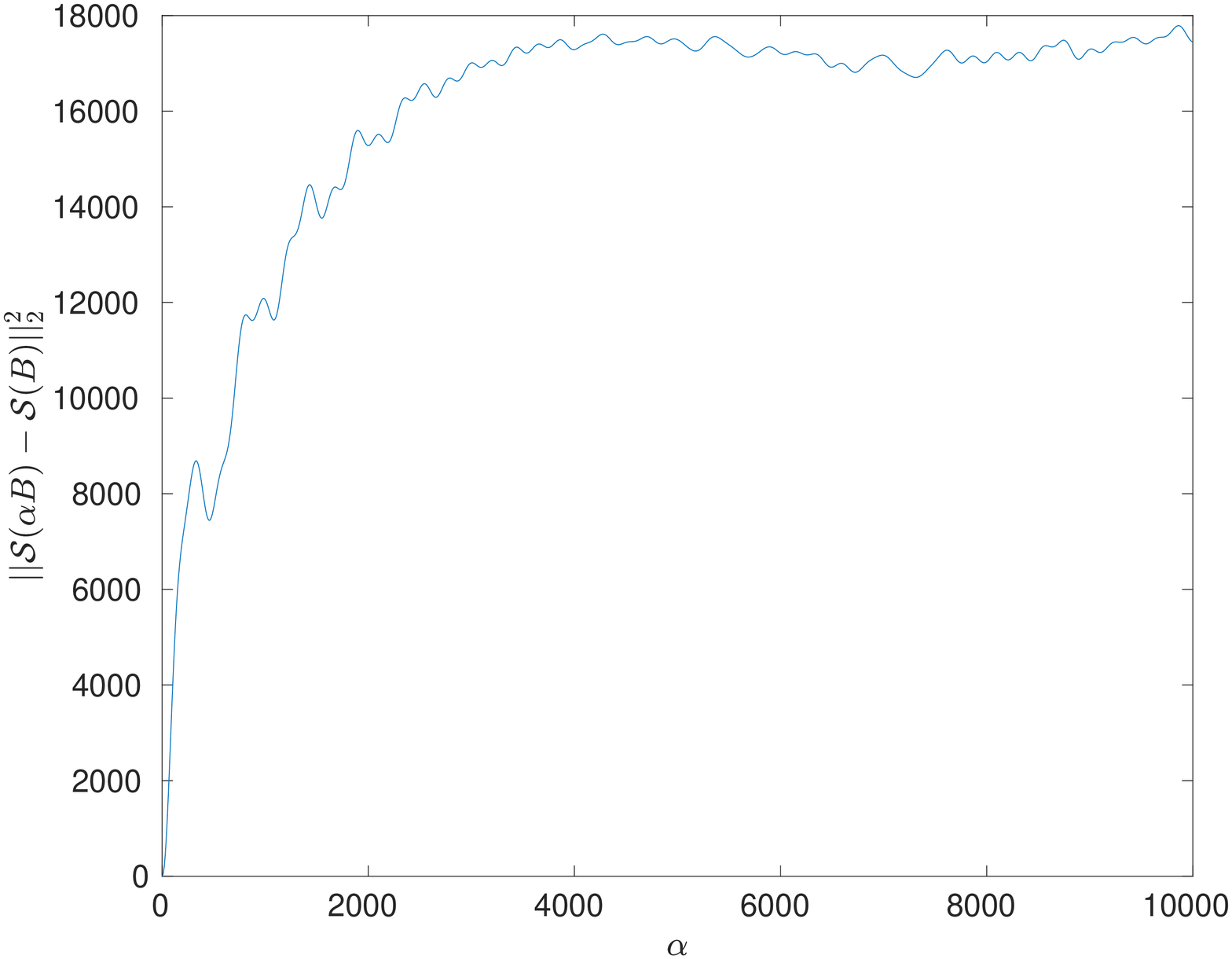}
\caption{\label{fig:non-convex} Plot of the residual error  $||\mathcal{S}(\alpha B)- \mathcal{S}(B) ||^2$ against the parameter $\alpha$ where $B$ is the solenoid example. Note that for sufficiently large $\alpha$ the curve is no longer convex.  } 
\end{figure}

We have already seen in Section \ref{section:ForwardProblemFormulation} that the Gateaux derivative of $\mathcal{S}$ can be calculated as a matrix attenuated ray transform. In a discretized setting this gives us an explicit way to calculate the Jacobian matrix of the ray data with respect to the values of $B$ in voxels. Preliminary experiments using regularized CGLS to solve the resulting linear system showed that for large magnetic fields the algorithm could converge to a local minimum, but not global minimum, of f  the residual function. 

\section{Conclusions}

We have an algorithm for  reconstructing magnetic fields, providing the field is weak enough, or the velocity great enough to not cause phase wrapping. Even though this is for simulated data, we have seen reconstruction is possible for experimental data, \cite{3DPNT}. We have suggested a general condition for which  
MNKM worked, namely the correct magnetic field can be retrieved if the angle by which the neutron spin precesses is under $\pi/2$.  We know from Theorem \ref{thm:PlaneBReconFormula} that the solution is unique with any amount of phase wrapping, however as the objective function is non-convex this presents a challenge for optimization based methods. In practice one would typically have a wide range of neutron speeds, and so the problem of phase wrapping might be avoided by using the faster neutrons to resolve the phase wrapping for rays encountering strong magnetic fields.

Ideally a robust reconstruction algorithm is required which can take care of phase wrapping issues when dealing with strong magnetic domains with a field as much as  $1 T.$  Ultimately for application purposes, e.g. quantum mechanical effects in superconductors and imaging electromagnetic devices, a  
further advancement is required. More often than not, the magnetic materials that we would like to image are composed of disjoint domains with a sharp change in the magnetic field at their boundary. To derive an innovative reconstruction algorithm for such a problem would really advance the field.

The uniqueness result of \cite{paternain2012attenuated} assumes smoothness of the metric, which is given in the Euclidean case. An important research topic for this applied problem is to reduce the smoothness assumption on the magnetic field (Higgs field in the geometric formulation). It is also an open problem to find an explicit non-linear reconstruction method.

There are various other methods that can be used to improve our technique. This problem could be formulated as a Bayesian inverse problem with assumed prior probability densities based on other physical knowledge. For example if domains are known from attenuation tomography, fields can be assumed to be close to constant within domains. 

\section*{Acknowledgements}
Parts of this work were supported by the European Union INTERREG \"{O}resund-Kattegat-Skagerrak fund, Otto M\o{}nsteds Fond, Royal Society Wolfson Research Merit Award and EPSRC grants EP/J010456/1, EP/M022498/1, EP/M010619/1. WL and ND would like to thank Gregory Eskin, Roman Novikov, Sean Holman and Gabriel Paternain for helpful discussions on non-Abelian ray transforms.

\section*{References}

\bibliography{thesis}

\bibliographystyle{plain}

\end{document}